\documentclass[final]{elsart}

\usepackage{amsmath}
\usepackage{amssymb,bm}
\usepackage{graphicx}
\usepackage{color}
\usepackage{multirow}
\usepackage{mathrsfs} 
\makeatletter
\newif\if@restonecol
\makeatother

\usepackage{algorithm}
\usepackage{algpseudocode}
\DeclareMathAlphabet\mathpzc{OT1}{pzc}{m}{it}
\let\mathcal=\mathpzc
\def\E{{\mathbb E}}
\def\P{{\mathbb P}}

\let\trueiiint=\iiint
\def\iiint{\mathop{\textstyle\trueiiint}\limits}
\def\intinfty{\int\limits_{\!\!-\infty\,\,}^{\,\,\infty\!\!}\kern-0.0em}
\def\iintinfty{\mathop{\int\!\!\int}\limits_{\!\!-\infty\,\,}^{\,\,\infty\!\!}\kern-0.0em}
\def\iiintinfty{\mathop{\int\!\!\int\!\!\int}\limits_{\!\!-\infty\,\,}^{\,\,\infty\!\!}\kern-0.0em}

\let\<=\langle
\let\>=\rangle
\let\^=\hat
\def\~#1{{\mbox{\sf#1}}}
\let\'=\vec
\let\-=\mathbf

\def\circ{\ifmmode\mathchar"220E\else$\mathchar"220E$\fi}
\def\@#1{{\cal #1}}

\def\COV{\mathrm{COV}}

\numberwithin{equation}{section}

\journal{Elsevier}

\begin{document}
\centerline{}
\begin{frontmatter}



\title{A surrogate accelerated multicanonical Monte Carlo method for uncertainty quantification}


\author[authorlabel1]{Keyi Wu}
\ead{wukeyi@sjtu.edu.cn}

\address[authorlabel1]{Department of Mathematics, Zhiyuan College, Shanghai Jiao Tong University,
Shanghai 200240, China.}

\author[authorlabel2]{Jinglai Li}
\ead{jinglaili@sjtu.edu.cn}
\address[authorlabel2]{Institute of Natural Sciences, Department of Mathematics, and MOE Key Laboratory of Scientific and Engineering Computing, Shanghai Jiao Tong University, Shanghai 200240, China. (Corresponding author)}


\medskip
\begin{center}
\end{center}

\begin{abstract}
In this work we consider a class of uncertainty quantification problems 
where the system performance or reliability is  characterized by a scalar parameter $y$. 
The performance parameter $y$ is random due to the 
presence of various sources of uncertainty in the system,  
and our goal is to estimate the probability density function (PDF) of $y$. 
We propose to use the multicanonical Monte Carlo (MMC) method,
 a special type of adaptive importance sampling algorithm,  to compute 
the PDF of interest.
Moreover, we develop an adaptive algorithm to construct local Gaussian process surrogates
to further accelerate the MMC iterations. 
With numerical examples we demonstrate that the proposed method can achieve
several orders of magnitudes of speedup over the standard Monte Carlo method. 
\end{abstract}


\end{frontmatter}

\section{Introduction}\label{s:intro}

Uncertainty is an inevitable feature of real-world engineering systems. In those systems uncertainty can rise from various of sources:
material properties, geometric parameters, boundary conditions, applied loadings and so on.
In practice, it is essentially important to characterize and quantify the impact of the uncertainties on the system performances, 
which constitutes a central task of the newly emerging field of Uncertainty Quantification (UQ).
To be specific, we consider the UQ problems in the following setting. 
We assume that the system is (formally) characterized by a performance function $y = g(\-x)$, where the input $\-x$ is a random vector collecting 
all the uncertain factors in the system and $y$ is a scalar indicating
the system performance or reliability (in what follows, we will simply refer to $y$ as the performance variable).  
A typical example is the structural design problems, in which $y$ can be the stress or the deformation. 
In this setting, the key task is to accurately assess and quantify the uncertainty in the performance parameter $y$. 
A challenge here is that real-world applications demand various statistical information of the performance $y$:
for example,  in robust design, the interests are mainly in the lower moments, especially the mean and the variance~\cite{du2004sequential},
in reliability analysis, it is mainly the tail probability~\cite{rackwitz2001reliability}, in risk management, one can be interested in the tail probability as well as some extreme quantiles~\cite{rockafellar2000optimization},     
and in utility optimization, the complete distribution of the performance parameter is required~\cite{hazelrigg1998framework}.  
To this end, a unified solution is to acquire the knowledge of the probability distribution of the performance parameter, which provides a complete characterization
of the uncertainty in it. 
In theory, the distribution of $y$ can be estimated by crucial Monte Carlo (MC) simulations, provided that a sufficient number of samples
can be afforded.  
In reality, however, the function $g: \-x\,\rightarrow\,y$ generally admits no analytical form, and evaluating function $g(\-x)$ must be done
by performing computer simulation of the underlying system, which renders 
 estimating the distribution of $y$ with crucial MC impractical. 

The main purpose of this work is to provide an efficient method to compute the full distribution of $y$. 
The proposed method has two major ingredients. 
First, we propose to sample the distribution of $y$ with the multicanonical Monte Carlo (MMC) method, which 
can be regarded as a more efficient alternative to MC. 
The MMC method was initially developed by Berg and Neuhaus~\cite{berg1991multicanonical,berg1992multicanonical} to explore the energy landscape of a given physical system, and later it has been adopted 
to simulate rare events, such as transmission errors in optical communication systems~\cite{holzloohner2003use,yevick2002multicanonical},
and the rare growth factors in random matrices~\cite{driscoll2007searching}. 
Roughly speaking, the MMC method constructs an iterative procedure that generates samples forming a flat histogram in the space of the parameter of interest (i.e., the energy in the original problem setup). 
As will be shown in Section~\ref{s:mmc}, the MMC method often requires to iterate many times and in each iteration it
employs Markov chain Monte Carlo (MCMC) simulations to draw a rather large number of samples. 
As a result, the direct use of MMC to sampling the distribution of the performance can still be computationally demanding, especially for systems with computationally intensive models. 
To this end, the second major component of our method is to employ computationally inexpensive surrogates to further reduce the computational cost of MMC.
In particular,  building on the method proposed in the work~\cite{conrad2014accelerating}, we adaptively construct local Gaussian process (GP) surrogates in the MCMC iteration.
{We choose to use the this method for the following reasons: 
first, the surrogate construction scheme is naturally incorporated in the MCMC iterations, which makes it convenient to use;
secondly, unlike many other surrogate based algorithms which introduce
 errors in the equilibrium distribution, this method samples asymptotically from the exact distribution of interest~\cite{conrad2014accelerating}. }

It should be noted that the purpose of the MMC method differs from that of the advanced sampling techniques developed in the field of reliability analysis or rare event simulations,  such as the cross entropy method~\cite{li2011efficient}, subset simulations~\cite{AuB99}, sequential Monte Carlo~\cite{cerou2012sequential}, etc. 
Namely, the purpose of those methods is to provide a variance-reduced estimator for 
for a specific parameter associated with the distribution of $y$, while that of our method is to estimate the distribution of $y$ itself.
{As will be shown in the next section, MMC is particularly useful for this purpose, which is our primary motivation to choose MMC over other advanced sampling schemes.}

The rest of this paper is organized as the following. 
We first review the MMC method in Section~\ref{s:mmc}, and then present our local GP construction algorithm in Section~3.
Finally numerical examples are provided in Section \ref{s:results} to demonstrate the
effectiveness of the proposed method.

\section{The multicanonical Monte Carlo method} \label{s:mmc}
In this section we introduce the MMC algorithm, largely following the presentation of \cite{bononi2009fresh}.
We start by summarizing the basic setup of our problem.
Let $\-x$ be a random vector taking values in the state space $X$, and 
 $y = g(\-x)$ be a real scalar function of $\-x$. 
 For simplicity we assume that both $\-x$ and $y$ are continuous random variables whose probability density functions exist. 
  We further assume that the PDF $p (\-x)$ of $\-x$ is known, possibly up to an unknown normalization constant,
 and our goal is to determine the PDF $\pi(y)$ of $y$.

\subsection{Flat histogram importance sampling}\label{sec:fh}
\begin{figure}
\centerline{\includegraphics[width=.9\textwidth]{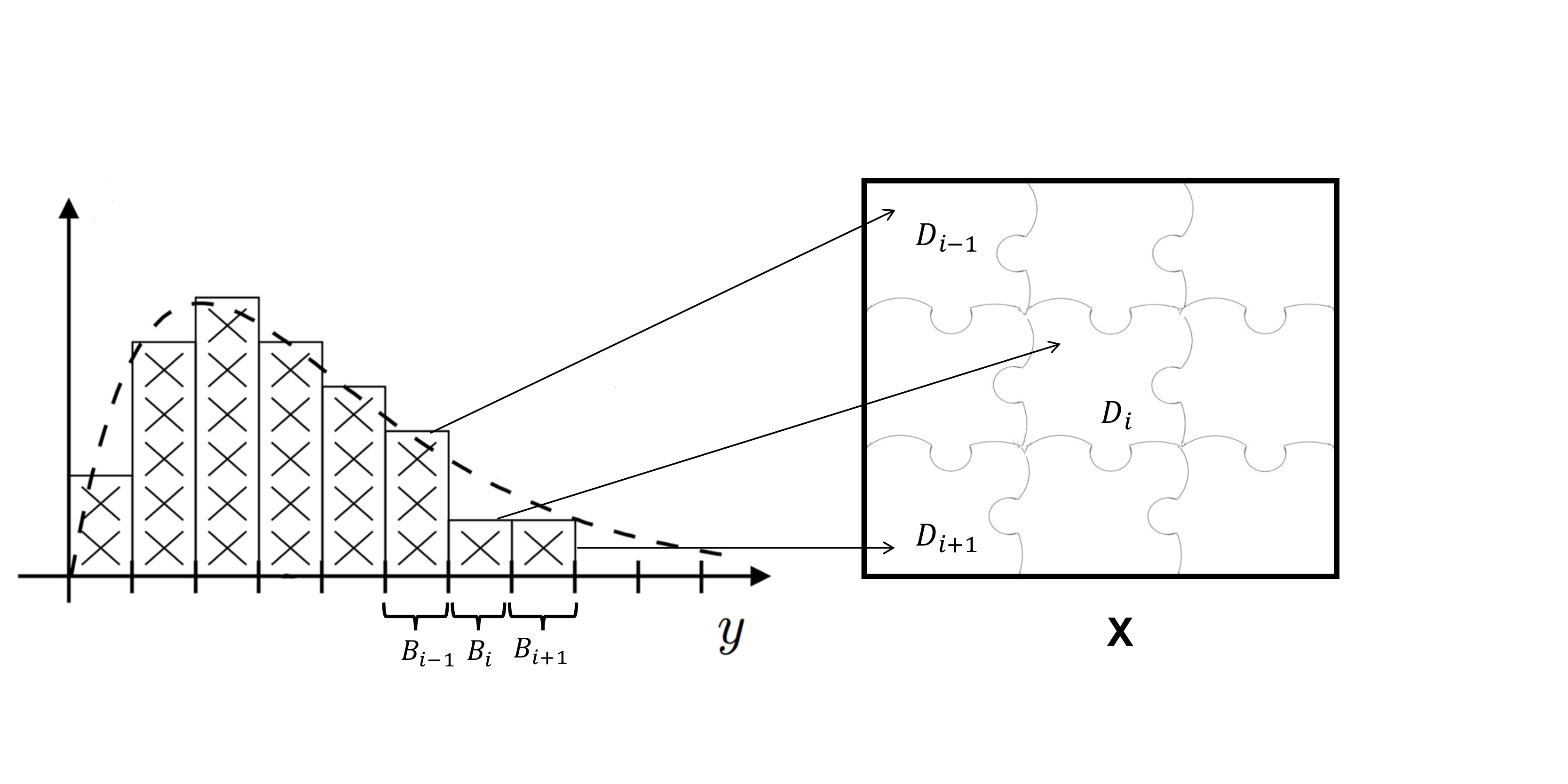}} \label{fig:bins}
\caption{Schematic illustration of the connection between $B_i$ and $D_i$.}
\end{figure}

A popular strategy to estimate the PDF of a continuous random variable $y$ with simulation,  
is to approximate the PDF with histograms, like a special case of the kernel density estimation. 
Suppose we are interested in the PDF of $y$ in a given closed interval $R_y$, and we first equally decompose $R_y$ into $M$ bins of width $\Delta$ centered at the discrete values $\{b_1, ... , b_M\}$. We define the $i$-th bin as the interval $B_i = [b_i - \frac{\Delta}{2}, b_i + \frac{\Delta}{2}]$ and the probability for $y$ to be in $B_i$ is $P_i =\P\{y \in B_i\}$.
The PDF of $y$ at point $y_i$ then can be approximated by \[\pi(y_i) \approx P_i / \Delta,\] if $\Delta$ is sufficiently small. This binning implicitly defines a partition of the input space $X$ into $M$ domains $\{ D_i \}^{M}_{i = 1}$, where
$$D_i = \{ \-x \in  X: g(\-x) \in B_i \}$$
is the domain in $X$ that maps into the $i$-th bin $B_i$. See Fig.~\ref{fig:bins} for an illustration. Note that, while $B_i$ are simple intervals, the domains $D_i$ are multidimensional regions with possibly tortuous topologies. 
As a result, the probability $P_i$ can be re-written as an integral in the input space:
\begin{equation}
P_i = \int_{D_i} p (\-x) dx = \int I_{D_i}(\-x) p (\-x) dx = \E[I_{D_i}(\-x)] ,\label{eq:w:1}
\end{equation}
where $I_{D_i}(\-x)$ is an indicator function defined as,
\begin{align*}
I_{D_i}(\-x) = \left\{
\begin{aligned}
 &1 ~~ \-x \in D_i;  \\
&0 ~~ \text{otherwise.} 
\end{aligned}
\right.
\end{align*}
Now suppose that $N$ samples $\{\-x_1,\ldots,\-x_N\}$ are drawn from the distribution $p(\-x)$, possibly with MCMC, 
$P_i$ can be evaluated with the MC estimator:
\begin{equation}
\hat{P}^{MC}_i = \frac{1}{N} \sum^{N}_{j = 1} I_{D_i}(\-x_j) = \frac{N_i}{N}, \label{eq:w:2}
\end{equation}
where $N_i$ is the number of samples that fall in bin $B_i$. 

As is well known, standard MC simulations have difficulty in reliably estimating the probabilities in the tail bins. 
The technique of importance sampling~(IS) can be used to address the issue.
Namely we choose a biasing distribution $q(\-x)$ and 
re-write \eqref{eq:w:1} as 
\begin{equation}
P_i = \int I_{D_i}(\-x) [\frac{p(\-x)}{q(\-x)}] q (\-x) d\-x = \E^{*} [I_{D_i}(X) w(X)] \label{eq:w:3}
\end{equation}
where  $w(\-x)  = p(\-x)/q(\-x)$ is the IS weight, and $\E^*$ indicates expectation with respect to the biasing distribution $q(\-x)$.
It follows that the IS estimator of $P_i$ becomes
\begin{equation}
  \hat{P}^{IS}_i  =  {\left( \frac{N^*_i}{N} \right)}  {\left[ \frac{1}{N^{*}_i} \sum^{N}_{j = 1} I_{D_i}(\-x_j)w(\-x_{j})\right]} \label{eq:w:5}
\end{equation}
where the samples $\{\-x_1,\ldots,\-x_N\}$ are now drawn from the biasing distribution $q(\-x)$,
and $N^*_i$ is the number of samples falling in region $D_i$.
For conciseness, we let $\hat{H}^*_i =   {\frac{N^*_i}{N} }$.
The intuition behind IS is that, the biasing distribution should assign higher probability in the region of interest 
than the original one, and thus it can draw more samples in that region. 

{ 
The key of IS is to choose an appropriate biasing distribution $q(\-x)$ that can help to achieve the objective of the simulation.
Unlike regular IS methods which usually employ biasing distributions that are easy to sample from, 
the MMC method chooses a biasing distribution $q(\-x)$  in the form of:
\begin{align}
 \label{eq:w:6}
q(\-x) = \left\{
\begin{aligned}
 & \frac{p(\-x)}{ \Theta(\-x)}   & \-x \in D;  \\
&0 & \-x\notin D. 
\end{aligned}
\right.
\end{align}
where $\Theta(\-x) = \Theta_i$. 
For $q(\-x)$ to be a well-defined distribution, we must have $\sum_{i=1}^M P_i/\Theta_i$=1.
It is easy to see that the distribution given in Eq.~\eqref{eq:w:6} assigns 
a constant weight to all $\-x \in D_i$: $w(\-x) = w_i$ for $\-x\in D_i$ where $w_i = \Theta_i$,  
which is referred to be as uniform-weight (UW).
In particular, if we let $\Theta_i= MP_i$ for all $x \in D_i, i = 1,...,M$.
the biasing distribution in Eq.~\eqref{eq:w:6} assigns equal probability to each bin
and zero probability for any region outside $D = \cup_{i=1}^M D_i$,  
namely, 
\begin{equation}
P^*_1=P_2^*=...P^*_M=1/M, 
\quad\mathrm{where}\quad
P^{*}_i = \int I_{D_i} (\-x) q(\-x) d\-x . \label{e:pis}
\end{equation}
We say such a biasing distribution as to be flat-histogram (FH).
FH is an important feature for our purpose which is to have a good estimate of $P_i$ for all $i=1\ldots M$.

 
 
\subsection{Multicanonical Monte Carlo}
It is easy to see, however, that the actual UW-FH distribution presented in Section~\ref{sec:fh} can not be used directly, as
$\Theta_i$ depend on the sought after unknown $P_i$. 
The MMC method addresses the issue in an incremental manner. Simply speaking MMC iteratively constructs a sequence of distributions 
\begin{align}
 \label{e:qn}
q_k(\-x) = \left\{
\begin{aligned}
 &\frac{p(\-x)}{\Theta_k(\-x)},    & \-x \in D;  \\
&0 & \-x\notin D, 
\end{aligned}
\right.
\end{align}
where $\Theta_k(\-x) = \Theta_{k,i}$ for $\-x \in D_i$, 
converging to the actual UW-FH distribution.
Specifically the sequence usually starts with $q_0(\-x)$ where $\Theta_{0,i} = \rho$ for all $i=1,\ldots,M$ and $\rho = \sum_{i=1}^M P_i\leq1$ is 
the probability that $y$ falls in the region of interest\footnote{ 
In practice, it is often convenient to assume that $\rho\approx1$ and in this case we have $q_0(\-x)\approx p(\-x)$.}. 
The iteration is then guided by the following equation: 
\begin{equation}
P^*_i = \int_{D_i} q(\-x)d\-x = \frac{\int_{D_i} p(\-x) d\-x}{c_\Theta \Theta_i} = \frac{P_i}{c_\Theta \Theta_i},   \label{e:theta_i2}
\end{equation}
or equivalently $P_i= P^*_i \Theta_i$.  
Namely, in the $k$-th iteration, one first draws $N$ samples $\{ \-x_j \} ^N_{j = 1}$ from the current distribution $q_k(\-x)$, 
and then updates $\{\Theta_{k+1,i}\}_{i=1}^M$ using the following formulas, which are derived from Eq.~\eqref{e:theta_i2},
\begin{subequations}
\label{e:params}
\begin{gather}
\hat{H}_{k,i} = \frac{N^*_{k,i}}{N},\\
P_{k,i} = \hat{H}_{k,i} * \Theta_{k,j},\\
\Theta_{k+1,i} = P_{k,i},
\end{gather}
\end{subequations}
where $N^*_{k,i}$ is the number of samples falling into region $D_i$ in the $k$-th iteration. 
We reinstate that, unlike a usual IS method, which often chooses a biasing distribution easy to sample from, the biasing distribution of the MMC method Eq.~\eqref{e:qn} is 
not a standard distribution, and thus directly sampling from the distribution is rather difficult. 
To this end, MMC usually employs MCMC algorithm to draw samples from $q_k(\-x)$. }
Formal convergence analysis, as well as possible improvements of the method are not discussed in this work,
and interested readers may consult, e.g., \cite{berg2000introduction,berg2004markov,iba2014multicanonical,landau2014guide},
and the references therein.

\section{Accelerating MMC with local GP surrogates} \label{s:surrogate}
In the MMC iteration, the main computational cost arises from performing the MCMC iteration to
draw samples from each $q_k(\-x)$, for each sample requires a full-scale simulation of the underlying system.
Thus, the MMC efficiency can be significantly improved by using computationally inexpensive surrogates in the MCMC scheme.
As is mentioned in Section~1, here we adopt the adaptive surrogate construction scheme developed in \cite{conrad2014accelerating}. 
In the work, the authors presented their method with two different surrogate models: the quadratic regression and the Gaussian processes, 
and their numerical results suggest that the GP model has better performance.   
We thus choose to use the GP model, while noting that other types of surrogates can also be used. 
In this section, we first briefly introduce the GP surrogate and then present the adaptive surrogate construction scheme modified for our specific use in MMC. 

\subsection{Gaussian process regression}

The GP surrogates, which are also known as kriging, have been widely used in many practical problems (see e.g.,~\cite{williams2006gaussian}).
The GP surrogate constructs the approximation of $g(\-x)$ in a nonparametric Bayesian regression framework~\cite{OHagan1978,williams2006gaussian}.
Specifically the target function $g(\-x)$ is cast as
\begin{equation}
g(\-x) = \mu_0(\-x)+\eta(\-x)
\end{equation}
 where $\mu_0(\-x)$ is a real-valued function and $\epsilon(\-x)$ is a zero mean Gaussian process whose covariance is specified by 
a kernel $K(\-x,\-x')$, namely, 
\[ \COV[\eta(\-x),\eta(\-x')] = K(\-x,\-x'). \] 
In practice, $\mu_0(\-x)$ can be represented as a linear or a quadratic polynomial whose coefficients can be determined by simple regression. 
In this work, we assume it is a quadratic polynomial.  
The kernel $K(\-x,\-x')$ is positive semidefinite and bounded. 
{Popular choices of the covariance functions include
squared exponential, exponential,  and Matern. 
The hyper-parameters inside the covariance functions can be prescribed or determined by maximizing the marginal likelihood function.  }
Suppose that $N$ computer simulations of the function $g(\-x)$ are performed 
 at parameter values $\-X^* := \left[\-x^*_1, \ldots \-x^*_n\right]$, yielding function evaluations $\-y^* := \left[ {y}^*_1, \ldots {y}^*_n\right]$, 
where \[{y}^*_i = g(\-x_i)\quad \mathrm{for} \quad i=1,\ldots,n.\] 
Suppose we want to predict the function values at a given point $\-x$, i.e.,
$y=g(\-x)$, the posterior of which is Gaussian:
\begin{equation}
  y ~|~\-x,\-X^*, \-y^* \sim\mathcal{N}(\mu(\-x), \sigma^2(\-x)). \label{e:post}
\end{equation}
{
\begin{subequations}
\label{e:gp}
The posterior mean of $y$ is 
\begin{equation}
\mu(\-x) =  \mu_0(\-x) + K(\-x,\-X^*)^T K(\-X^*, \-X^*)^{-1}(\-y^*-\mu_0(\-X^*)),\label{eq:postMean}
\end{equation}
and the posterior variance is 
\begin{equation}
\label{eq:postcovsimple} 
\sigma^2(\-x)= K(\-x,\-x) - K(\-x,\-X^*)^T \, K(\-X^*, \-X^*)^{-1} \, K(\-X^*,\-x), 
\end{equation}
\end{subequations}
where the notation $K(\-A,\-B)$ to denote the matrix of the covariance evaluated at all pairs of points in set $\-A$ and in set $\-B$~\cite{williams2006gaussian}.}
Eq.~\eqref{eq:postMean} can be used as the surrogate to predict the function values at points of interest, 
and Eq.~\eqref{eq:postcovsimple} provides a measure of confidence in the predicted values.

 \subsection{Local GP construction}
 In the standard GP methods, the surrogates are constructed with all the data points.  
Constructing the GP surrogate can be very costly when the data set becomes large,
as it involves inverting a large covariance matrix. On the other hand, it has been well noted that data points far from the point of interest have little  influence on the prediction (assuming the usual choices of covariance). 
Thus, a natural choice is to construct GP only with the data points near the point of interest. 
The resulting surrogate is thus \emph{local}, in the sense that it is only intended to be accurate at the point of interest. 
Next we discuss in detail how to construct a local GP surrogate at point $\-x$ given a collection of model evaluations: $\-S:=\{(\-x_i,y_i)\}_{i=1}^{n_S}$
where $y_i = g(\-x_i)$ for $i=1...n_S$.

First we need to determine how many data points we want to use in the surrogate construction.
Following the suggestion of \cite{conrad2014accelerating}, we choose the number of data points $n$ as
 \[n = \sqrt{d_x}(d_x + 1)(d_x + 2) / 2,\]
 where $d_x$ is the dimensionality of $\-x$. This choice allows us to have sufficient data points to perform a quadratic regression for $\mu_0(\-x)$.
 The specific points used to build the surrogate are chosen with the nearest neighbor  (NN) method: namely, we use the $n$ points closest to $\-x$ to construct the GP surrogate.
 It has been pointed out that the NN method only provides a suboptimal point selection, and better selection strategy can be obtained 
 by solving an optimization problem. 
 However, in our problem, the GP construction must be done repeatedly in the MCMC scheme, and as result
 even very fast optimization may significantly increase the total computational cost. 
 In this respect, we nevertheless adopt the NN method for the sake of computational {simplicity}. 
 In what follows, we refer to a local GP surrogate constructed with the prescribed procedure, as $\widetilde{g}(\-x |\mathbf{S})$.

 \subsection{MCMC with local GP surrogates} \label{sec:algo}

In this section, we present a modified version of the  local surrogate accelerated MCMC scheme developed in \cite{conrad2014accelerating}. 
The method embeds an adaptive surrogate construction in the MCMC iteration: in each iteration the method constructs
a local surrogate using data set $\-S$, for the proposed point and the current point, and decides whether it needs model refinement; 
when refinement is needed, the algorithm then refines the surrogate by evaluating more points near the proposed point or the current one 
depending on where the refinement is triggered;
all the evaluated points are included in the data set $\-S$ which will be used for constructing surrogates in the next step. 
In \cite{conrad2014accelerating}, refinement is triggered by either of two criteria. The first is random: with probability $\gamma_t$, the
model refined at either the current point or the proposed point. We follow 
the random criterion in our algorithm as it is the essential for the theoretical convergence of the algorithm. 
The second criterion used in \cite{conrad2014accelerating}, intended to make the algorithm efficient in practice, 
is based on an error indicator of the acceptance probability. {In this work, we follow the random criteria and choose $\gamma_t$ to be a constant for simplicity.}
{
We use, however, a different practical criterion, taking advantage of the special structure of the target distribution $q_k(\-x)$ in Eq.~\eqref{e:qn}.
Namely, it is easy to see that, for $q_k(\-x)$ in Eq.~\eqref{e:qn}, an error in the surrogate does not 
cause an error in the acceptance probability unless the surrogate assigns the sample into a wrong bin, assuming a symmetric proposal distribution.}
{Specifically, suppose the current sample is $\-x^-$ and the proposed sample is $\-x^+$, 
and the posterior mean and variance of the GP at $\-x^+$ are $y^+$ and $\epsilon^2$ respectively.  
Suppose it is assigned to bin $B_i=[b_{i}-\Delta/2,\,b_i+\Delta/2]$ based on the predicted value $y^+$, 
and the probability that the assignment of $\-x_i$ is incorrect can be computed as
\begin{multline}
\beta(\-x^+):=\P[ g(\-x^+) <b_{i}-\Delta/2\, \mathrm{or}\, g(\-x^+)>b_i+\Delta/2]\\ = \Phi(b_{i}-\Delta/2,y^+,\epsilon)-\Phi(b_{i}+\Delta/2,y^+,\epsilon)+1,
\end{multline}
where $\Phi(\cdot, y^+,\epsilon)$ is the cumulative density function (CDF) of a normal distribution with mean $y^+$ and 
standard deviation $\epsilon$.
Thus we can define the refinement criteria as that the misassignment probability $\beta$ is smaller than a threshold value: 
$\beta <\beta_{\max}$. {Since the refinement criteria is applied to each iteration, 
the probability that the acceptance probability computed with the surrogate is erroneous is bounded by $2\beta_{\max}$,
in any iteration. }
{As a result, to achieve this probability boundedness, we only need to check if $\-x^+$ satisfies the quality condition: $\beta(\-x^+)<\beta_{\max}$, as $\-x^-$ has been verified in the previous iteration. }
We outline our algorithm in Algorithm~1, where the surrogate construction is
integrated into a standard Metropolis-Hastings (MH) MCMC scheme. 


\begin{algorithm}
  \label{alg:mh}
    \caption{Metropolis-Hastings with local GP surrogates}
    \label{alg:thm4}
    \begin{algorithmic}[1]
	\For {$t = 1,...,T$}
	\State $(\mathbf{x}_{t+1},y_{t+1},\mathbf{S}_{t+1}) \gets K_t(\mathbf{x}_t, y_t, \mathbf{S}_t, q(\cdot;\-y_t),\gamma_t,a_{\min})$
	\EndFor
       \\
        \Procedure{$K_t$}{$\mathbf{x}^{-}, y^-, \mathbf{S},  q(\cdot;y^-),\gamma$, $a_{\min}$}
	\State Draw proposal $\bf{x}^{+} \sim {\Pi}(\mathbf{x}^{-},\cdot )$
	\State $(y^{+} ,\epsilon^+)\gets \widetilde{g}(\mathbf{x}^{+}, \mathbf{S})$ \label{line}
	\If {$u \sim \mathrm{Uniform}(0,1) < \gamma$}
	       \State $y^+=g(\mathbf{x}^+)$
		\State $\mathbf{S} \gets \mathbf{S}\cup\{(\mathbf{x}^+,y^+)\}$
	\Else
	\State $\beta\gets 1 +\Phi(y_1,y^+,\epsilon)-\Phi(y_2,y^+,\epsilon^+)$
	\If {$\beta>\beta_{\max}$}
			       \State $y^+=g(\-x^+)$
		\State $\mathbf{S} \gets \mathbf{S}\cup\{(\mathbf{x}^+,y^+)\}$	
	\EndIf
	\EndIf
                \State $\alpha\gets q(\mathbf{x}^+;y^-)/q(\mathbf{x}^-;y^-)$
	 \If {$u \sim \mathrm{Uniform}(0,1) <\alpha$} 
	\State\Return $(\mathbf{x}^{+}, y^+,\mathbf{S})$
	\Else 
	\State \Return $(\mathbf{x}^{-}, y^-,\mathbf{S})$
	\EndIf

        \EndProcedure       
       
   \end{algorithmic}
\end{algorithm}

We have the following remarks regarding the proposed algorithm,
highlighting its differences from that given in \cite{conrad2014accelerating} in addition to the refinement criteria. 
\begin{itemize}
\item
As a pre-processing of the first MMC iteration, we choose $n_o$ points, and use them as the initial data set $\-S$. 
These points can be chosen in many different ways: sampling according to $p(\-x)$, Latin hypercube, or experimental design methods. 
For the succeeding MMC iterations, the data set $\-S$ is simply taken to be that obtained in the previous round. 
\item Unlike regular MCMC methods, in each iteration our algorithm returns the sample $\-x_t$ as well as the function value $y_t$ for the sample. 
Note that,  the function values are needed in Eqs.~\eqref{e:params}, 
and thus by recording the function values, we can compute Eq.~\eqref{e:params} without evaluating the function again.
\item 
{As has been discussed in the beginning of Section~\ref{sec:algo}}, in each iteration we only need to consider the quality of the surrogate at the proposed point $\-x^+$ thanks to the special
structure of $q_k(\-x)$, while in the original algorithm,
 both $\-x^+$ and $\-x^-$ need to be examined. 
\item In our algorithm, when model refinement is needed, we simply evaluate the current point $\-x^+$. It has been suggested 
that this strategy may lead to poor conditioned regression in particular when polynomial surrogates are used,
and as an alternative a space filling approach is used in \cite{conrad2014accelerating}.
However, we have found it is not a very serious issue for the GP surrogates in our numerical tests,   
and, considering that the space filling method requires an extra optimization step, we choose to directly evaluate $\-x^+$ for {simplicity}'s sake. 
\end{itemize}





{Finally we note that it is an very interesting problem to analyze the convergence property of the algorithm. 
To this end, the convergence analysis in \cite{conrad2014accelerating} can provide certain useful results of the MCMC iterations.
However, since the algorithm is a combination of the two methods, 
a formal convergence analysis can be very challenging, and so is not pursued in this work.}

\section{Numerical examples} \label{s:results}
We use three numerical examples to demonstrate the performance of the proposed GP accelerated MMC (GP-MMC) method. 
Before proceeding to the examples, we describe the specific GP surrogate used in all the three examples.
First in all the examples we use an anisotropic covariance function in the form of:
\begin{equation}
K(\-x,\-x') = a \exp\left[ -\sum_{i=1}^{d_x}\frac{|x_i-x'_i|^p}{l_i}\right], \label{e:corrfun}
\end{equation}
where $p$ is a prescribed positive integer which usually takes values of $1$ (the exponential kernel) or $2$ (the squared exponential kernel), the coefficient $a$ is determined with empirical Bayes in the iteration,
 and the correlation length $\-l = (l_1,...,l_{d_x})$ is determined from the initial data set
 and is not adjusted in the iteration. 
Note that, the correlation length $\-l$ can also be determined with empirical Bayes in the iteration if desired, 
but we choose not to do so here for simplicity's sake,  as it requires to numerically solving 
an optimization problem. 
  
\subsection{A multi-dimensional analytical example}

Our first example is a multi-dimensional problem where the performance function is 
\[g(\-x) = \min\{g_1(\-x), \, g_2(\-x)\}-1,\]
with 
{\[g_1(\-x) = \left \lVert  \-x - \-x_1 \right \rVert,\quad\mathrm{and}\quad g_2(\-x) = \left \lVert  \-x - \-x_2 \right \rVert.\]
The input $\-x$ are multidimensional independently distributed standard normal random variables  and $\-x_1,\-x_2$ are two fixed points}. 
It is obvious that each $D_i$ has two possibly disjoint sections: $\{\-x\,|\, g_1(\-x)\in B_i\}$
and $\{\-x\,|\, g_2(\-x)\in B_i\}$, which makes the problem challenging for many variance-reducing sampling techniques. 

{We first test our method for the two dimensional case and choose $\-x_1=(3,3)$ and $\-x_2=(3,-3)$ respectively.}
We run standard MC simulations with $10^{7}$ samples, and use its results as the ``truth''
 to validate the estimates of the MMC methods.
In the first numerical experiment, we perform MMC simulations without using surrogates, where $10$ iterations are used with $10^5$ samples in each iteration,
resulting in a total computational cost of $10^6$ full-model simulations. 
When constructing the PDF, we use $R_y=[-1, 54]$ which is divided into $55$ bins.
In Fig.~\ref{fig:Hs} we show the histograms obtained in the 1st, 2nd and the final MMC iteration, from which one can see that the histograms
tend to become flat as the iterations proceed. 
\begin{figure}
\centerline{\includegraphics[width=.33\textwidth]{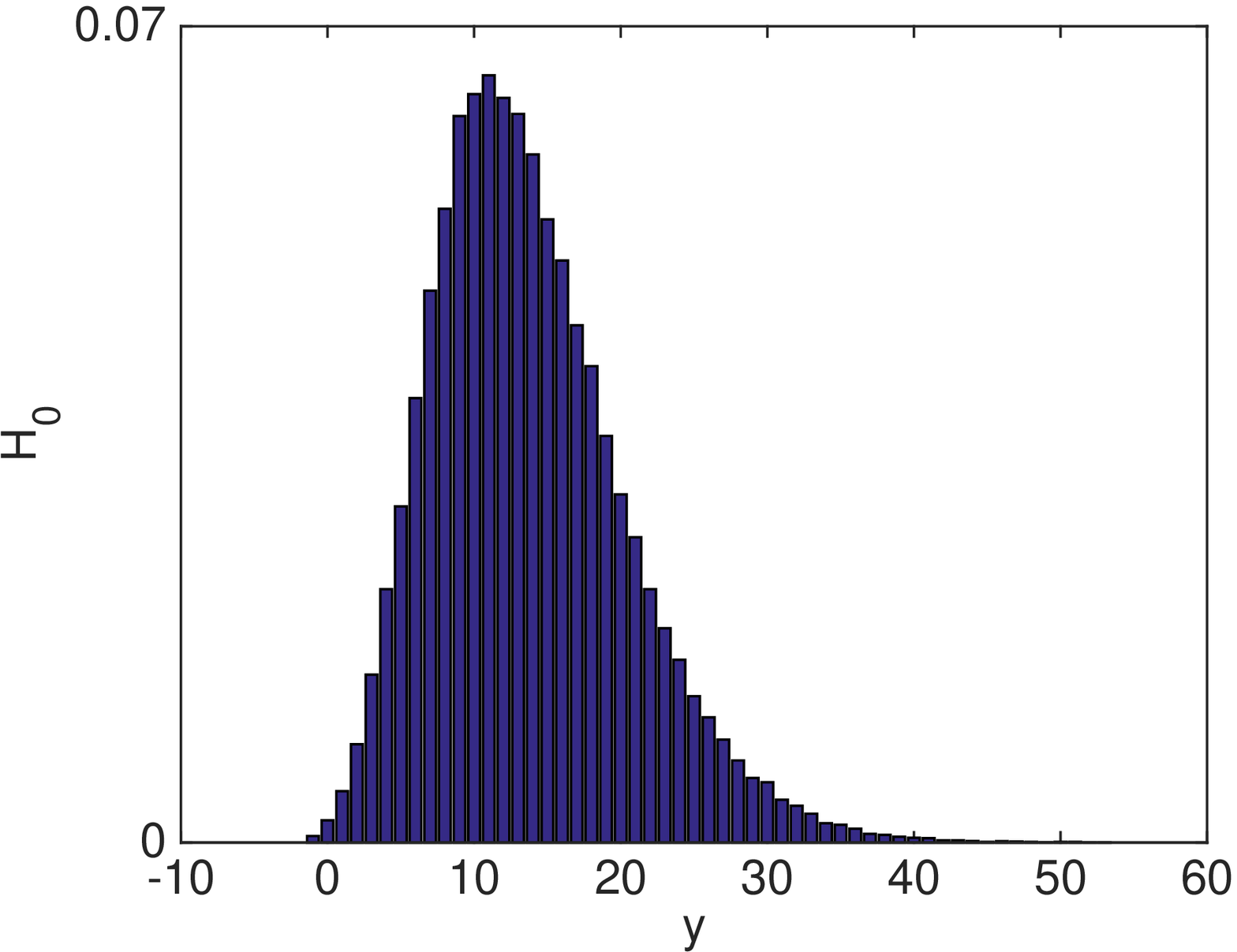}
\includegraphics[width=.33\textwidth]{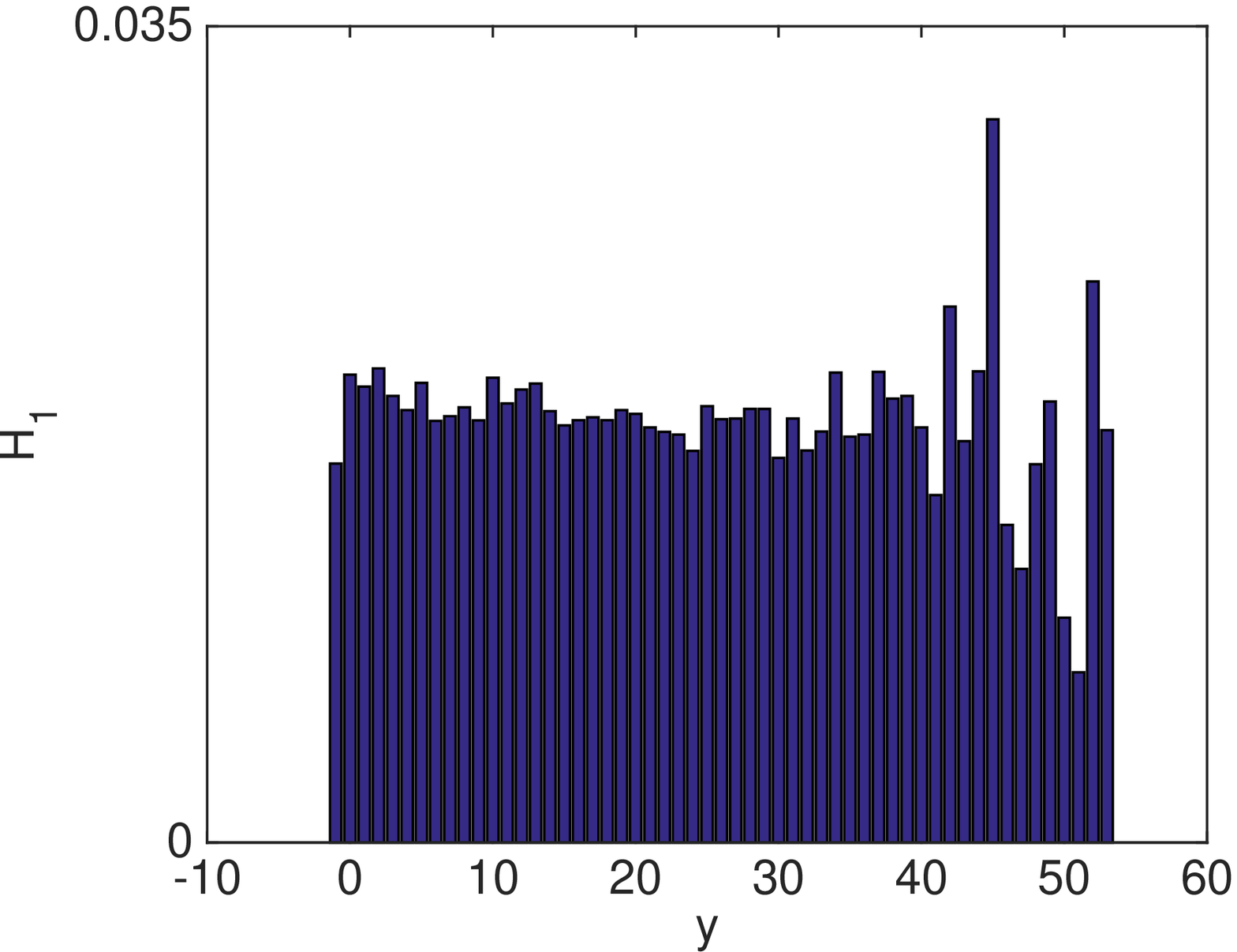}
\includegraphics[width=.33\textwidth]{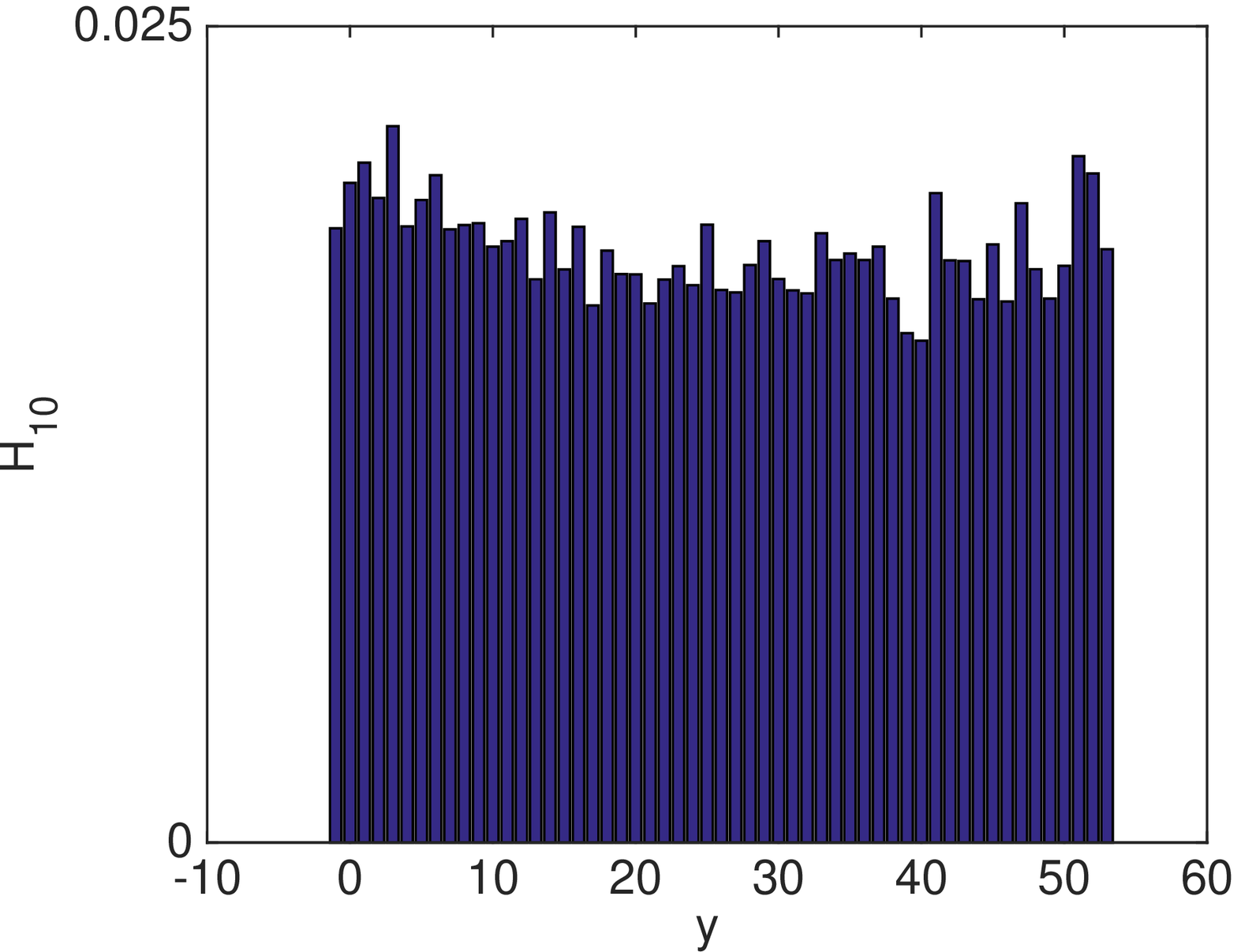}} \label{fig:Hs}
\caption{The histograms of the first two steps and the 10th iteration of MMC.}
\end{figure}

Our second numerical experiment is to run MMC with the assistance of the GP surrogates, and, as is in the first experiment, 
we again use $10$ iterations with $10^5$ samples in each.  
In the GP-MMC computation, we construct the GP surrogates as is described in the beginning of the section, where the kernel is given by Eq.~\eqref{e:corrfun} with $p=1$. 
 The initial data set contains 50 samples randomly drawn from the distribution of $\-x$, 
and we choose the random model refinement probability $\gamma_t =10^{-4}$.
The key parameter in the algorithm is the maximum misassignment probability $\beta_{\max}$, 
and to examine the robustness of our method against the choices of $\beta_{\max}$, we implement our method with various values of $\beta_{\max}$ and show the results in Table~\ref{ta:results1}.
\begin{table}[htbp]
\centering
\begin{tabular}{lcccccc}
\hline
$\beta_{\max}$  & 0.92&0.76&0.32&0.05&0.003 &plain MMC\\
\hline
true model evals&796&810&809&926&1089&$10^6$\\
\hline
maximum {RelErr}&0.1775& 
0.148&
0.1058&
0.1321&
0.1217&0.0921
\\
\hline
average {RelErr}& 0.0419&
0.039&
0.0327&
0.0333&
0.0345&0.0225
\\
\hline
\end{tabular} 
\caption{(Example 1) The performance results of GP-MMC with various values of $\beta_{\max}$.}
\label{ta:results1}
\end{table}
In particular, for the results of each value of $\beta_{\max}$, we show the number of true model evaluations, the maximum  
and the average relative errors (compared to the MC results) of all the bins.

\begin{figure}
\centerline{\includegraphics[width=.85\textwidth]{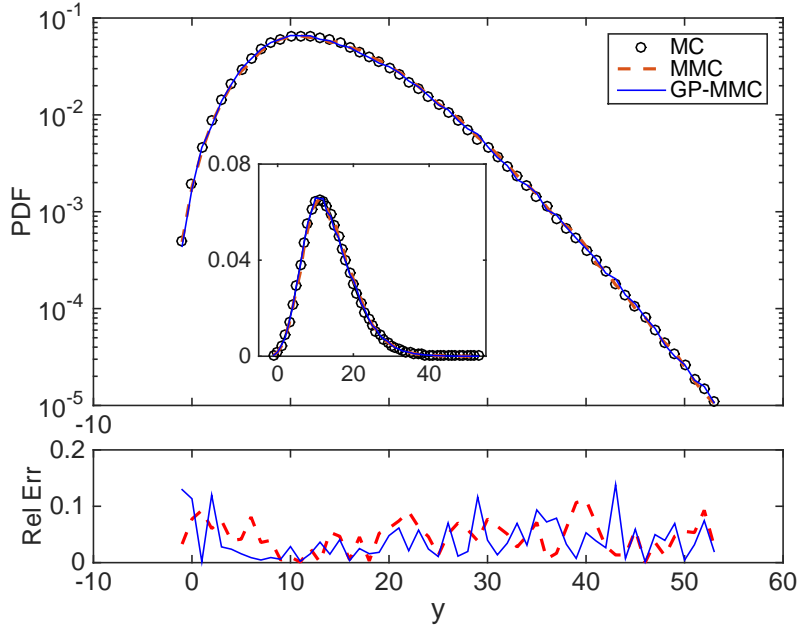}} 
\caption{(Example 1) Top: the PDF of $y$ obtained by MC (circles), MMC (dashed line) and GP-MMC (solid line) on a logarithmic scale;
inset is the same plots on a linear scale. 
Bottom: the relative error in the PDF obtained by MMC (dashed) and GP-MMC (solid).}\label{fig:pdf}
\end{figure}


One can see that, the method performs well even for a very large misassignment probability, 
and the results are rather robust for different values of $\beta_{\max}$ except that the number of true model evaluations 
grows as $\beta_{\max}$ becomes smaller. 

To further compare the results, we plot the PDF obtained by MC, MMC and GP-MMC with $\beta_{\max}=0.05$,  in Fig.~\ref{fig:pdf} (Top),
and one can see that the results of the three methods agree very well with each other. 
To have a quantitative assessment of the performance, we compute the relative error of the MMC and the GP-MMC estimates,
against the results of plain MC:
\begin{equation}
\mathrm{RelErr}_\mathrm{MMC} =\frac{|\hat{p}_\mathrm{MMC} -\hat{p}_\mathrm{MC} |}{\hat{p}_\mathrm{MC}},\quad
\mathrm{RelErr}_\mathrm{GPMMC} =\frac{|\hat{p}_\mathrm{GPMMC} -\hat{p}_\mathrm{MC} |}{\hat{p}_\mathrm{MC}},
\end{equation}
and show the results in Fig.~\ref{fig:pdf} (Bottom).

We see that, the relative errors in both MMC and GP-MMC are around $0.1$, indicating 
that both MMC and GP-MMC produce reliable estimates of the PDF of $y$.
To further compare the performance, we computed the mean, the variance, the 3rd, the 4th and the 5th central moments of $y$
using the samples obtained by the three methods
shown in Table~\ref{ta:moments_eg1}, which shows that the results obtained by the three methods agree well with each other.
Regarding the computational cost, the MMC method uses $10^6$ full model evaluations while our GP-MMC method
only uses  less than a thousand full-model evaluations. 

\begin{table}[htbp]

\centering
{
\begin{tabular}{lcccccc}
\hline
moment  & mean&var&3rd&4th&5th\\
\hline
MC&14.21&43.58&217.42&7340.55&108583.52\\
\hline
MMC&14.43&44.04 &241.11&7505.02&113171.57\\
\hline
GP-MMC& 14.28&44.04& 230.10&7456.33&111877.63\\
\hline
\end{tabular} 
\medskip

\caption{(Example 1) The mean, variance, and 3rd--5th central moments of $y$, estimated
by MC, MMC and GP-MMC .}
\label{ta:moments_eg1}}
\end{table}

{
We also consider the performance of the proposed method with respect to different sample sizes and dimensionality.
To this end, we first perform the GP-MMC method as well as standard MMC with different number of samples in each iteration,
in which $\beta_{\max}$ is taken to be $0.075$.  
The results are shown in Table~\ref{ta:eg1_size}, and as expected, with more samples in each iteration, the results become 
more accurate at the price of more true model evaluations. 
Next, we consider the example with different number of dimensions. 
In this case we let $\-x_1=(1,...,1)^d$ and $\-x_2 = (-1,...,-1)^d$ for $d=2,\,8,\,16$.
We perform both MMC and GP-MMC in each case, and in the GP-MMC 
 we take $\beta_{\max} =0.075$ and the number of sample size in each iteration to be $5\times10^4$. 
 The performance comparison is shown in Table~\ref{ta:eg1_dim}.
We can see from the results that as the dimensionality increases, the GP-MMC method requires more true model evaluations,
but the computational cost saving compared to standard MMC is still significant even for the case of 16 dimensions.  
 Overall we have found that the performance of the GP-MMC method is rather robust with respect to the sample size and the dimensionality. 
 }
\renewcommand{\multirowsetup}{\centering}  
\begin{table}[htbp]
{
\centering
\begin{tabular}{c|ccccc}
\hline
 & sample size  & 1e+4&1e+5&1e+6\\

\hline
 \multirow{2}[5]{2cm}{MMC}
&maximum RelErr&0.3097&0.1858&0.0466\\\cline{2-5}
&average RelErr&0.0791&0.0387 &0.0120\\
\hline

 \multirow{3}[6]{2cm}{GP-MMC}
&true model evals&891&1855&2033\\\cline{2-5} 
&maximum RelErr&0.2593&0.0906&0.0576\\\cline{2-5}
&average RelErr&0.087&0.0328 &0.0147\\
\hline
\end{tabular} 
\medskip

\caption{(Example 1) The performance of MMC and GP-MMC with respect to various sample sizes.}
\label{ta:eg1_size}
}
\end{table}

\renewcommand{\multirowsetup}{\centering}  

\begin{table}[htbp]

{
\centering
\begin{tabular}{c|ccccc}
\hline
&dimension  & 2&8&16\\
\hline

 \multirow{2}[5]{2cm}{MMC}
&maximum RelErr&0.1173&0.1168 &0.1531\\\cline{2-5} 
&average RelErr&0.0225&0.0370 &0.0566\\

\hline

 \multirow{2}[6]{2cm}{GP-MMC}
 
&true model evals&891&3226&16886\\\cline{2-5} 
&maximum RelErr&0.1497&0.1521&0.1692\\\cline{2-5} 
&average RelErr&0.0414&0.0422&0.0665\\

\hline
\end{tabular} 
\medskip

\caption{(Example 1) The performance of MMC and GP-MMC with respect to various numbers of dimensions.}
\label{ta:eg1_dim}
}
\end{table}

\subsection{Cantilever beam}
We now consider a cantilever beam problem~\cite{li2011efficient,wu1990advanced}
 as illustrated in
Fig.~\ref{fig:beam}, with width $w$, height $t$, length $L$,
and subject to transverse load $Y$ and horizontal load $X$.
This is a popular benchmark problem in the reliability analysis literature, where 
the performance function is
$$
y =\frac{4L^3}{Ewt}\sqrt{\left(\frac{Y}{t^2}\right)^2+\left(\frac{X}{w^2}\right)^2},
$$
which represents the deflection of the beam. 
In this example, we assume that  the beam length is fixed $L=100$, and the beam width $w$, the height $x$, the applied loads $X$ and $Y$, 
as well as the elastic module $E$ of the material, are random parameters.  We further assume that these random parameters
are all independently distributed, with each following a normal distribution. 
The means and the variances of the parameters are summarized in Table~\ref{ta:params}. 
\begin{table}[htbp]
\centering
\begin{tabular}{lcccccc}
\hline
parameter  & w&t&X&Y&E&\\
\hline
mean&4&4&500&1000&$2.9\times10^{6}$&\\
\hline
variance& 0.001&   0.0001 & 100 &   100&   $1.45\times10^6$\\
\hline
\end{tabular} 
\medskip

\caption{(Example 2) The mean and variance of the random parameters in the cantilever beam model.}
\label{ta:params}
\end{table}

\begin{figure}
\centerline{\includegraphics[width=.65\textwidth]{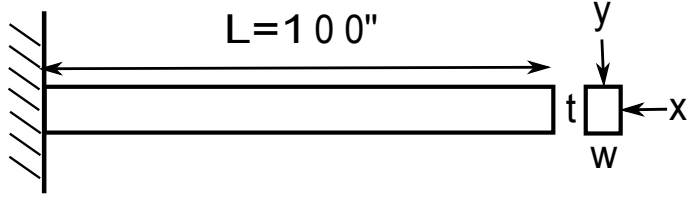}} \label{fig:beam}
\caption{(Example 2) Schematic illustration of a cantilever beam
  subject to horizontal and vertical loads.}
\end{figure}

In this example, we also compute the PDF of $y$ with three methods: plain MC, MMC and GP-MMC.
In the MC simulations, we use $10^9$ full model evaluations. 
In both MMC and GP-MMC, we use $10$ iterations where $10^5$ samples in each iteration.
In the GP-MMC computation, 
the number of initial data and the values of $\gamma_t$ are the same as those used in the first example.
The GP kernel is also given by Eq.~\eqref{e:corrfun} with $p=1$. 
Also, we test the GP-MMC method with various values of $\beta_{\max}$ and show the results in Table~\ref{ta:results_beam}.
In this example, we use $R_y=[0.56, 0.66]$ which is divided into $40$ bins.
To compare the results, we plot the PDF obtained by MC, MMC and GP-MMC with $\beta_{\max}=0.32$ which requires $4775$ true model evaluations,  
as well as the relative errors of MMC and GP-MMC, in Figs.~\ref{fig:compare_beam}.
We also show the same moment plots as is in the first example in Table~\ref{ta:moments_beam_data}.
All the figures indicate that our GP-MMC method yields very reliable estimates of the PDF of $y$,
while its computational cost is significantly lower than both MC and standard MMC.

\begin{table}[htbp]

\centering
\begin{tabular}{lcccccc}
\hline
$\beta_{\max}$  & 0.92&0.76&0.32&0.05&0.003\\
\hline
true model evals&894&2523&4775&7456&7589&\\
\hline
maximum {RelErr}&0.116&0.143 &0.102&0.099&0.089\\
\hline
average {RelErr}& 0.034&0.042& 0.038&0.038&0.037\\
\hline
\end{tabular} 
\medskip

\caption{(Example 2) The performance results of GP-MMC with various values of $\beta_{\max}$.}
\label{ta:results_beam}
\end{table}

\begin{figure}
\centerline{\includegraphics[width=.85\textwidth]{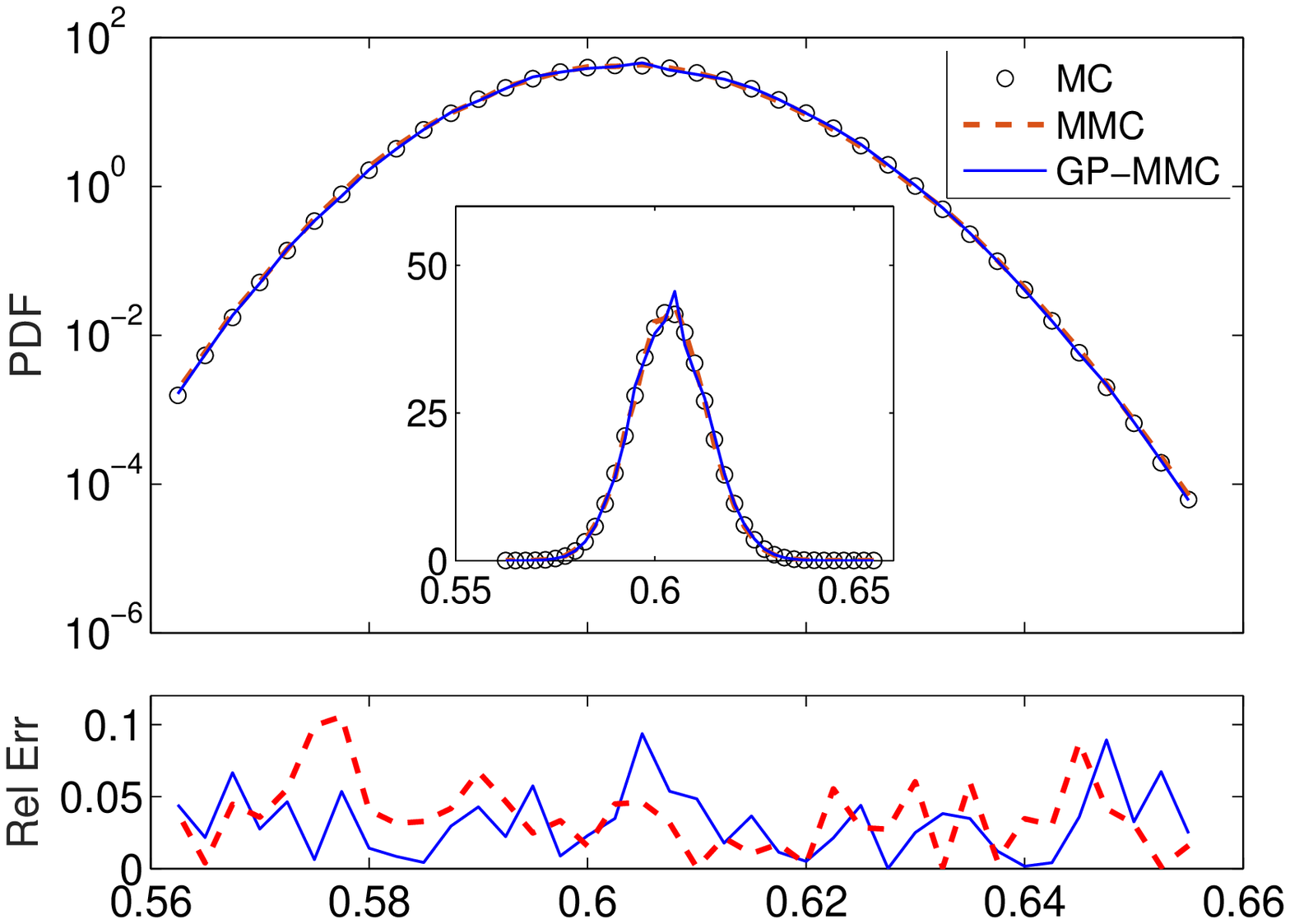}} 
\caption{(Example 2) Top: the PDF of $y$ obtained by MC (circles), MMC (dashed line) and GP-MMC (solid line) on a logarithmic scale;
inset is the same plots on a linear scale. 
Bottom: the relative error in the PDF obtained by MMC (dashed) and GP-MMC (solid).}\label{fig:compare_beam}
\end{figure}

\begin{table}[htbp]

\centering
{
\begin{tabular}{lcccccc}
\hline
moment  & mean&var&3rd&4th&5th\\
\hline
MC&0.6024&8.99e-5&6.28e-8&2.43e-8&4.96e-11\\
\hline
MMC&0.6024&8.97e-5 &3.04e-8&2.43e-8&3.32e-11\\
\hline
GP-MMC& 0.6025&9.04e-5& 7.55e-8&2.46e-8&5.54e-11\\
\hline
\end{tabular} 
\medskip

\caption{(Example 2) The mean, variance, and 3rd--5th central moments of $y$, estimated
by MC, MMC and GP-MMC .}

\label{ta:moments_beam_data}}
\end{table}

\subsection{Random PDE example}
Finally we consider a random partial differential equation (PDE) example: 
 a  two-dimensional Poisson equation on region $\Gamma = [0,1]\times[0,1]$:  
 \begin{subequations}
\label{e:poisson}
 \begin{gather}
 \nabla(a(\-x)\nabla u(\-x)) = f(\-x),\\
u = 0 \quad \mathrm{on} \quad \partial\Gamma,
\end{gather}
\end{subequations}
where $a(\-x)$ is a random field and $\partial\Gamma$ is the boundary of $\Gamma$. 
We want to compute the statistical distribution of the value of $u$ at location $\-x_*\in \Gamma$. 
 A physical interpretation of the problem is the following:  we consider a steady flow in an isotropic aquifer subject to random permeability~\cite{anderson1992applied}, 
and we are interested in the statistical information of the hydraulic head at a particular location $\-x_*$.

\begin{figure}
\centerline{\includegraphics[width=.6\textwidth]{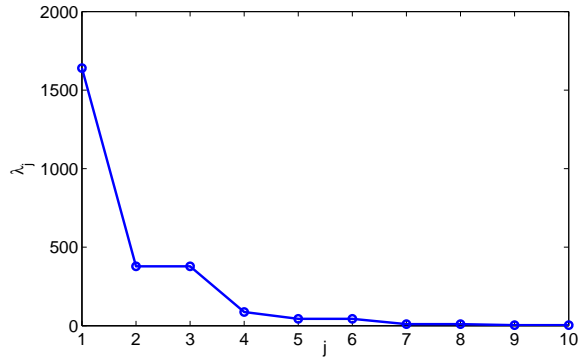}} 
\caption{(Example 3) The eigenvalues of the KL expansion plotted in a descending order.}\label{fig:kl_eig}
\end{figure}

\begin{figure}
\centerline{\includegraphics[width=.49\textwidth]{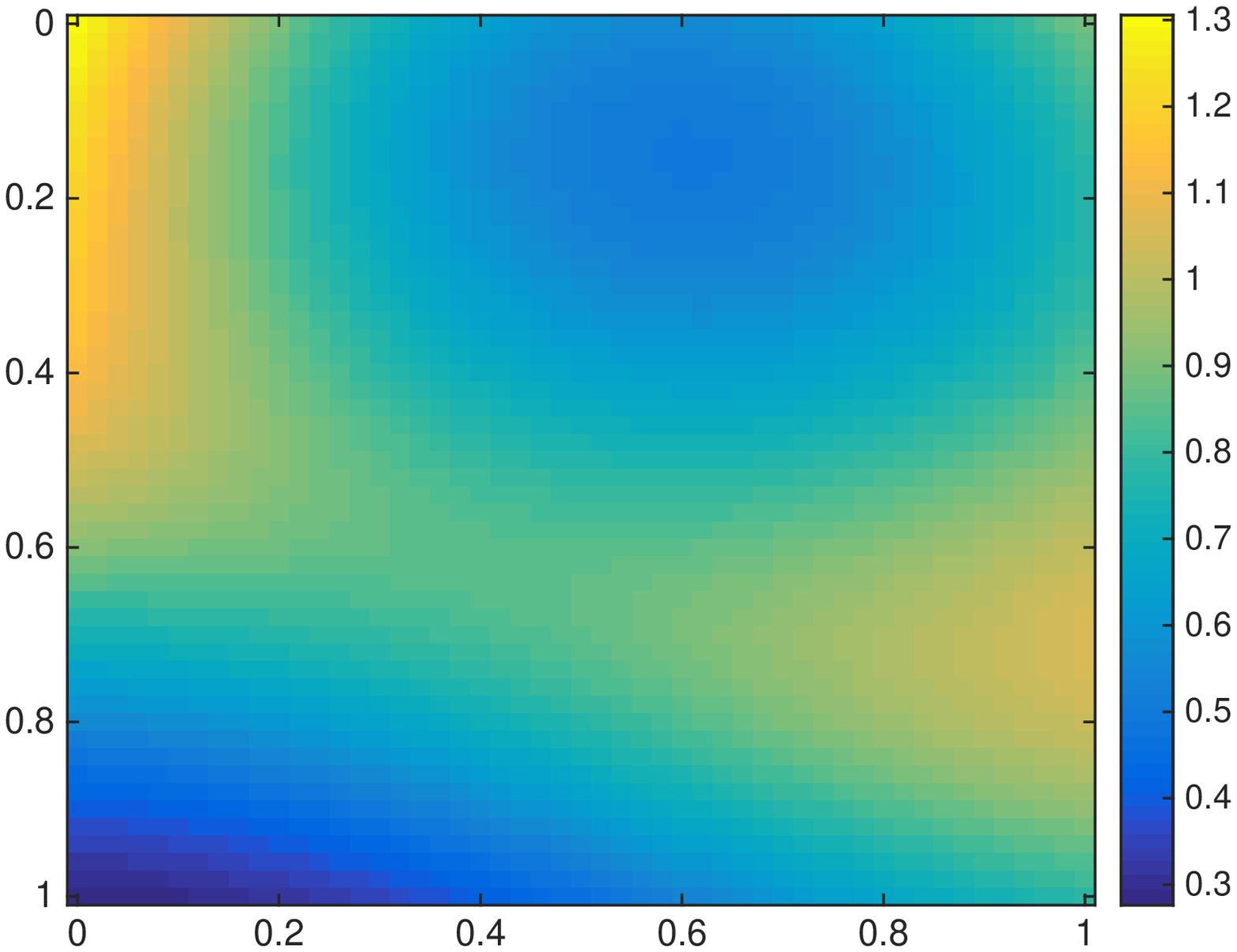}
\includegraphics[width=.49\textwidth]{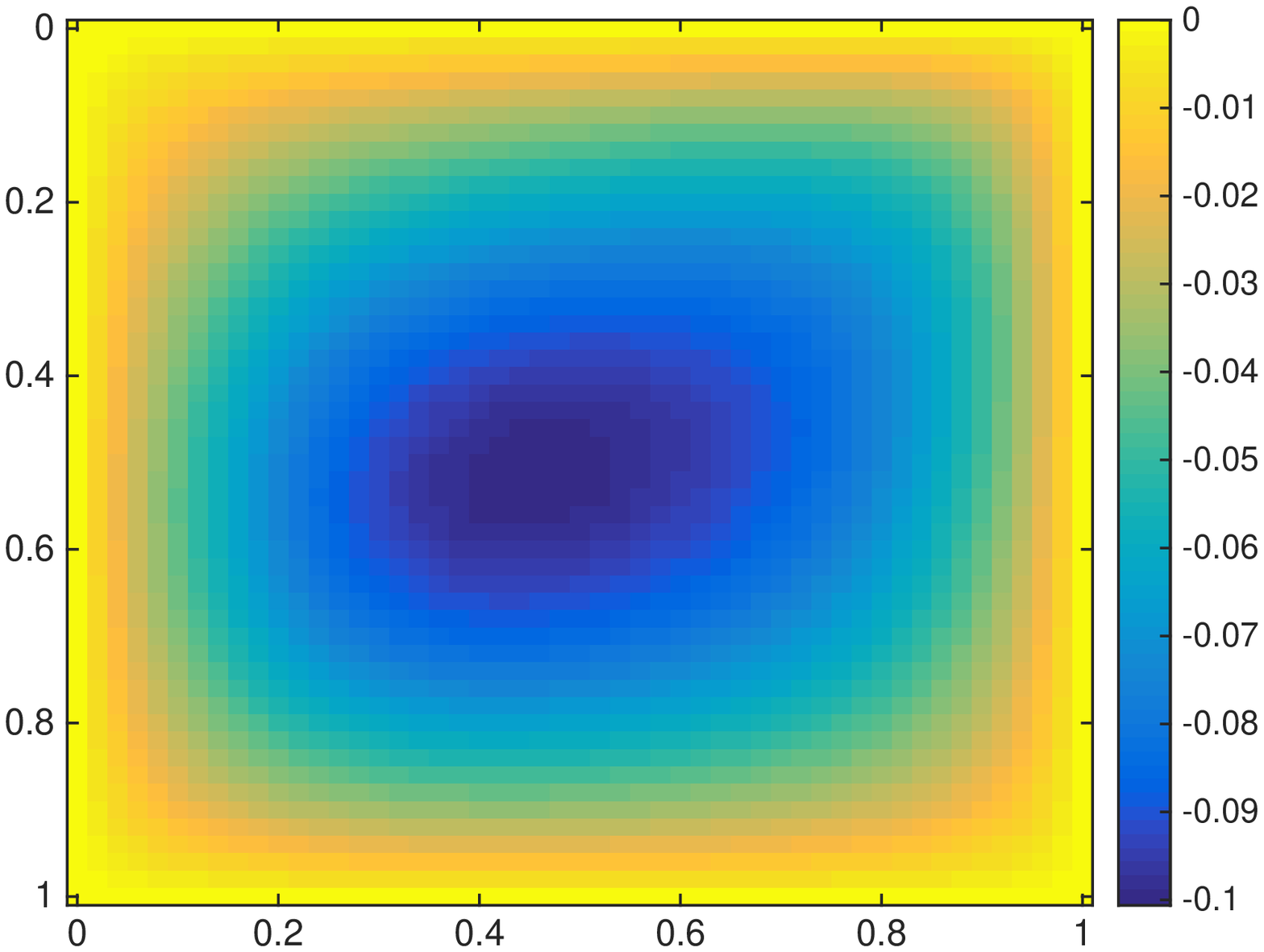}} \label{fig:samp_pde}
\caption{Left: a randomly drawn coefficient sample $a(\-x)$. Right: the solution of Eq.~\eqref{e:poisson} associated with $a(\-x)$.}
\end{figure}

We further assume the permeability is a log-normal random field, namely, $a(\-x) = a_o\exp(z(\-x))$ where $z(\-x)$ is a Gaussian random field with zero mean and covariance kernel,
\begin{equation}
\Sigma(\-x_1,\-x_2)= \exp(-\frac{\|\-x_1-\-x_2\|^2}{\Delta}).
\end{equation}
In this example we take $a_0=1$ and $\Delta =0.6$.
In practice, the random field $z(\-x)$ in the PDE is often represented with a truncated Karhunen-Lo\`eve~(K-L) expansion. 
Namely, let $\{\lambda_j,\xi_j(\-x)\}_{j=1}^\infty$ be the eigenvalue-eigenfunction pairs of the covariance kernel $\Sigma(\cdot,\cdot)$ such that
$\lambda_j>\lambda_{j+1}$ for all $j=1...\infty$, 
and we can approximate $z(\-x)$ with 
\begin{equation}
z(\-x) = \sum_{j=1}^J c_j\sqrt{\lambda_j} \xi_j(\-x), \label{e:kle}
\end{equation}
where $\-c = (c_1,...,c_J)$ follows a standard isotropic normal distribution.
Thus the dimensionality of the problem is reduced to $J$ and in this example we choose $J=10$.
We plot the eigenvalues associated with the 10 KL modes in a descending order in Fig.~\ref{fig:kl_eig},
which suggests that 10 KL-modes can sufficiently represent the Gaussian field $z(\-x)$ in this problem.
Moreover, in the numerical simulations, we take $f(\-x)=1$ and $\-x^*=(0.5,\,0.5)$.
A sample coefficient $a(\-x)$ and the associated solution $u(\-x)$ is shown in Fig.~\ref{fig:samp_pde}.

\begin{figure}
\centerline{\includegraphics[width=.85\textwidth]{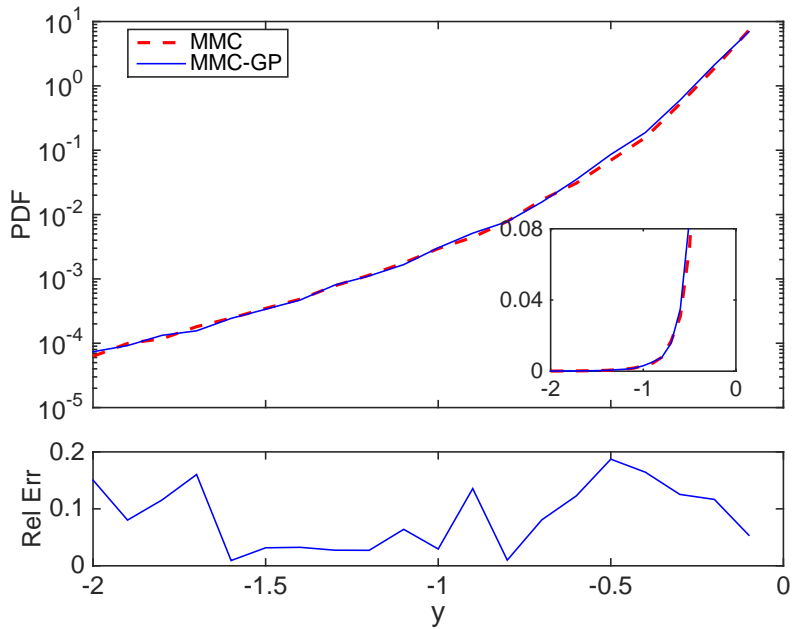}} 
\caption{(Example 3) Top: the PDF of $y$ obtained by MMC (dashed line) and GP-MMC (solid line) on a logarithmic scale;
inset is the same plots on a linear scale. 
Bottom: the relative error in the PDF.}\label{fig:compare_pde}
\end{figure}

As the computational cost for solving Eq.~\eqref{e:poisson} is rather high, which renders standard MC unfeasible, 
 we choose to only perform MMC and GP-MMC simulations in this problem. 
In both cases, we use 10 iterations with 20000 samples in each iteration. 
In GP-MMC, we use the covariance function~\eqref{e:corrfun} with $p=2$.
The number of initial samples is 400, $\gamma_t=10^{-4}$ and $\beta_{\max} = 0.05$. 
As a result the total number of true model evaluations is 4885.
When constructing the PDF, we use $R_y=[-2, 0]$ divided into $20$ bins.
We plot the PDF computed with MMC and GP-MMC as well as the relative error in the two results in Fig.~\ref{fig:compare_pde}.
One can see from the figures that the results of GP-MMC agree very well with those of plain MMC, while it only uses around one fortieth 
true model evaluations of the plain MMC.  

\section{Conclusions}
We consider a special type of UQ problems where the system performance is characterized by a scalar parameter. 
We propose to use a MMC based method to compute the distribution of the performance parameter,
and we also propose to use a local GP surrogate to accelerate the MMC simulations.  
Based on the work~\cite{conrad2014accelerating}, we design an adaptive algorithm that can effectively refine the GP surrogate in the MMC iterations. 
With numerical examples, we demonstrate that the proposed GP-MMC method can efficiently and accurately 
compute the distribution of the performance parameter.  
We expect the proposed method can be useful in various fields of applications, such as reliability analysis, risk management, and utility optimizations. 

There are a number of possible improvements and extensions of the proposed method that we plan to investigate in the future. 
First there are some well-known open issues with GP: most notably, how to choose
the best covariance functions, and such a choice may certainly affect the performance of our MMC-GP method. 
To this end, we hope to develop approaches that can effectively choose the covariance functions for our MMC method. 
{Second, as has been mentioned in Section~3, we are not able to provide a convergence analysis of the proposed method in this paper and 
we hope to address the issue in a future work. }
Third we are also interested in more general uncertainty propagation problems where the output is a multidimensional vector rather than a scalar. 
In this case, the standard MMC scheme does not apply directly, due to the multi-dimensionality of the output.
We plan to tackle such problems with modified MMC algorithms. 
Finally, we note that the Wang-Landau algorithms, which can be regarded as a variant of MMC, have been 
applied to the Bayesian inference problems (e.g. \cite{chopin2012free}), and we hope that our GP-MMC method can be applied to such problems as well.
In this case, we expect that our method can further improve the computational efficiency of the Bayesian inferences,
 thanks to the use of surrogates.

\section*{Acknowledgment}
We want to thank the anonymous reviewer for his/her very constructive comments.
The work was partially supported by the National Natural Science Foundation of China under grant number 11301337.

\bibliographystyle{plain}
\bibliography{mmc}

\end{document}